\documentclass[12pt,a4paper]{amsart}
\usepackage{amsfonts}
\usepackage{amsthm}
\usepackage{amsmath}
\usepackage{amscd}
\usepackage[latin2]{inputenc}
\usepackage{t1enc}
\usepackage[mathscr]{eucal}
\usepackage{indentfirst}
\usepackage{graphicx}
\usepackage{graphics}
\usepackage{pict2e}
\usepackage{epic}
\usepackage{url}
\usepackage{enumerate}
\numberwithin{equation}{section}
\usepackage[margin=2.9cm]{geometry}

\usepackage{tikz}
\usetikzlibrary{shapes}
\tikzstyle{vertex}=[draw=black,circle,fill=black,minimum size=4pt, inner sep=0pt, outer sep=0pt,text=white,line width=0mm]
\tikzstyle{c0}=[shape=circle, minimum size=4pt, fill=white]
\tikzstyle{c1}=[shape=rectangle, minimum size=7pt, fill=red]
\tikzstyle{c2}=[shape=diamond, minimum size=10pt, fill=blue]

\theoremstyle{plain}
\newtheorem{Th}{Theorem}[section]

\newtheorem{Cor}[Th]{Corollary}

 \theoremstyle{definition}
\newtheorem{Def}[Th]{Definition}
\newtheorem{Conj}[Th]{Conjecture}
\newtheorem{Rem}[Th]{Remark}
\newtheorem{?}[Th]{Problem}

\renewcommand{\tilde}{\widetilde}
\newcommand{\HWR}{H_{\mathrm{WR}}}
\newcommand{\Hind}{H_{\mathrm{ind}}}
\newcommand{\V}[1]{V_{#1}}
\newcommand{\E}[1]{E_{#1}}
\newcommand{\A}[1]{A_{#1}}
\newcommand{\B}[1]{B_{#1}}

\begin{document}

\title[Widom-Rowlinson model and regular graphs]{The Widom-Rowlinson model, the hard-core model and the extremality of the complete graph}

\author[E.~Cohen]{Emma Cohen}
\address{School of Mathematics, Georgia Institute of Technology\\ Atlanta, GA 30332-0160}
\email{ecohen32@math.gatech.edu}

\author[P.~Csikv\'ari]{P\'{e}ter Csikv\'{a}ri}
\address{Massachusetts Institute of Technology \\ Department of Mathematics \\
Cambridge MA 02139 \&  MTA-ELTE Geometric and Algebraic Combinatorics Research Group
\\ H-1117 Budapest
\\ P\'{a}zm\'{a}ny P\'{e}ter s\'{e}t\'{a}ny 1/C \\ Hungary} 
\email{peter.csikvari@gmail.com}

\author[W.~Perkins]{Will Perkins}
\address{School of Mathematics, University of Birmingham, UK}
\email{william.perkins@gmail.com}

\author[P.~Tetali]{Prasad Tetali}
\address{School of Mathematics and School of Computer Science, Georgia Institute of Technology\\ Atlanta, GA 30332-0160}
\email{tetali@math.gatech.edu}

\thanks{The second author  is partially supported by the National Science Foundation under grant no. DMS-1500219, by the MTA R\'enyi "Lend\"ulet" Groups and Graphs Research Group, by the ERC Consolidator Grant 648017, and by the Hungarian National Research, Development and Innovation Office, NKFIH grant K109684. Research of the last author is supported in part by the NSF grant DMS-1407657.}

 \subjclass[2010]{Primary: 05C35. Secondary: 05C31, 05C70, 05C80}

 \keywords{graph homomorphisms, Widom-Rowlinson model, hard-core model} 

\begin{abstract} Let $\HWR$ be the path on $3$ vertices with a loop at each vertex. D.~Galvin \cite{Gal1,Gal2} conjectured, and E.~Cohen, W.~Perkins and P.~Tetali \cite{CPT} proved that for any $d$-regular simple graph $G$ on $n$ vertices we have
$$\hom(G,\HWR)\leq \hom(K_{d+1},\HWR)^{n/(d+1)}.$$
In this paper we give a short proof of this theorem together with the proof of a conjecture of Cohen, Perkins and Tetali \cite{CPT}. Our main tool is a simple bijection between the Widom-Rowlinson model and the hard-core model on another graph. We also give a large class of graphs $H$ for which we have
$$\hom(G,H)\leq \hom(K_{d+1},H)^{n/(d+1)}.$$
In particular, we show that the above inequality holds if $H$ is a path or a cycle of even length at least $6$ with loops at every vertex.
\end{abstract}

\maketitle

\section{Introduction} For graphs $G$ and $H$, with vertex and edge sets $\V{G}, \E{G}, \V{H}$, and $\E{H}$ respectively, a map $\varphi:\V{G}\to \V{H}$ is a homomorphism if $(\varphi(u),\varphi(v))\in \E{H}$ whenever $(u,v)\in \E{G}$. The number of homomorphisms from $G$ to $H$ is denoted by $\hom(G,H)$. When $H=\Hind$, an edge with a loop at one end, homomorphisms from $G$ to $\Hind$ correspond to independent sets in the graph $G$, and so $\hom(G,\Hind)$ counts the number of independent sets in $G$.

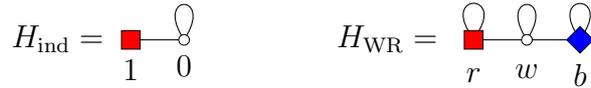
\begin{figure}[h]
\centering
\[\Hind = 
\begin{tikzpicture}[scale=.7, baseline=-1mm]
  \node[vertex, c1, label={below:$1$}] (v1) at (-1,0) {};
  \node[vertex, c0, label={below:$0$}] (v0) at (0,0) {};
  \draw (v1) to (v0) .. controls +(60:1) and +(120:1) .. (v0);
\end{tikzpicture}
\qquad
\qquad
\HWR = 
\begin{tikzpicture}[scale=.7, baseline=-1mm]
  \node[vertex, c1, label={below:$\vphantom{b}r$}] (v1) at (-1,0) {};
  \node[vertex, c0, label={below:$\vphantom{b}w$}] (v0) at (0,0) {};

  \node[vertex, c2, label={below:$b$}] (v2) at (1,0) {};
  \draw (v1) .. controls +(60:1) and +(120:1) .. (v1) to
        (v0) .. controls +(60:1) and +(120:1) .. (v0) to
        (v2) .. controls +(60:1) and +(120:1) .. (v2);
\end{tikzpicture}
\]
\caption{The target graphs for the Widom-Rowlinson model and the hard-core model.}
\label{fig:H-examples}
\end{figure}

For a given $H$, the set of homomorphisms from $G$ to $H$ correspond to valid configurations in a corresponding statistical physics model with \emph{hard constraints} (forbidden local configurations).  The independent sets of $G$ are the valid configurations of the \emph{hard-core model} on $G$, a model of a random independent set from a graph. Another notable case is when $H=\HWR$, a path on $3$ vertices with a loop at each vertex. In this case, we can imagine a homomorphism from $G$ to $\HWR$ as a $3$-coloring of the vertex set of $G$ subject to the requirement that a blue and a red vertex cannot be adjacent (with white vertices considered unoccupied); such a coloring is called a \emph{Widom-Rowlinson configuration} of $G$, from the Widom-Rowlinson model of two particle types which repulse each other \cite{WR,BHW}. See Figure \ref{fig:H-examples}.

For a fixed graph $H$, it is natural to study the normalized graph parameter 
$$p_H(G):=\hom(G,H)^{1/|\V{G}|},$$
where $\V{G}$ denotes the number of vertices of the graph $G$.

For $H=\Hind$, J.~Kahn \cite{Kahn} proved that for any $d$-regular bipartite graph $G$, 
$$p_{\Hind}(G)\leq p_{\Hind}(K_{d,d}),$$
where $K_{d,d}$ is the complete bipartite graph with classes of size $d$. Y.~Zhao \cite{Zhao1} showed that one could drop the condition of bipartiteness in Kahn's theorem. That is, he showed that $p_{\Hind}(G)\leq p_{\Hind}(K_{d,d})$, for \emph{any} $d$-regular graph $G$. Y.~Zhao proved his result by reducing the general case to the bipartite case with a clever trick. He proved that
$$p_{\Hind}(G)\leq p_{\Hind}(G\times K_2),$$
where $G\times K_2$ is the bipartite graph obtained by replacing every vertex $u$ of $\V{G}$ by a pair of vertices $(u,0)$ and $(u,1)$ and replacing every edge $(u,v)\in \E{G}$ by the pair of edges $((u,0),(v,1))$ and  $((u,1),(v,0))$. This is clearly a bipartite graph, and if $G$ is $d$-regular then $G\times K_2$ is still $d$-regular. 

D.~Galvin \cite{Gal1,Gal2} conjectured a different behavior for $H = \HWR$: that instead of $K_{d,d}$, the complete graph $K_{d+1}$ maximizes $p_{\Hind}(G)$ among $d$-regular graphs $G$. E.~Cohen, W.~Perkins and P.~Tetali \cite{CPT} proved that this was indeed the case:

\begin{Th} \label{WR} \cite{CPT} For any $d$-regular simple graph $G$ on $n$ vertices we have
$$p_{\HWR}(G)\leq p_{\HWR}(K_{d+1});$$
in other words,
$$\hom(G,\HWR)\leq \hom(K_{d+1},\HWR)^{n/(d+1)}.$$
\end{Th}

One of the goals of this paper is to give a very simple proof of this fact\footnote{In fact, Theorem~\ref{WR} follows from a stronger result in \cite{CPT} that the Widom-Rowlinson {\em occupancy fraction} is maximized by $K_{d+1}$.  We note that this stronger result also follows from the transformation below and Theorem 1 of \cite{davies2015independent}.}, along with a slight generalization. We use a trick similar to that used by Y.~Zhao \cite{Zhao1,Zhao2}. We will need the following definition:

\begin{Def} The \emph{extended line graph} $\tilde{H}$ of a (bipartite) graph $H$ has $\V{\tilde{H}} = \E{H}$; two edges $e$ and $f$ of $H$ are adjacent in $\tilde{H}$ if 
\begin{enumerate}[(a)]
  \item $e=f$,
  \item $e$ and $f$ share a common vertex, or
  \item $e$ and $f$ are opposite edges of a $4$-cycle in $G$.
\end{enumerate}
\end{Def} 
Throughout, $\V{H}$ and $\E{H}$ refer to the vertex-set and edge-set, respectively, of the graph $H$. If $H$ is bipartite, we use $\A{H}$ and $\B{H}$ to refer to the parts of a fixed bipartition.
Now we can give a generalization of Theorem~\ref{WR}:

\begin{Th} \label{gen}  If $\tilde{H}$ is the extended line graph of a bipartite graph $H$, then for any $d$-regular simple graph $G$ on $n$ vertices we have
$$p_{\tilde{H}}(G)\leq p_{\tilde{H}}(K_{d+1}),$$
or in other words,
$$\hom(G,\tilde{H})\leq \hom(K_{d+1},\tilde{H})^{n/(d+1)}.$$
\end{Th}

To see that Theorem~\ref{gen} is a generalization of Theorem~\ref{WR} it suffices to check that $\HWR$ is precisely the extended line graph of the path on $4$ vertices.
In Section~\ref{extension} we will prove a slight generalization of Theorem~\ref{gen} which allows for weights on the vertices of $H$.

\section{Short proof of Theorem~\ref{WR}}

We are not the first to notice the following connection between the Widom-Rowlinson model and the hardcore model (see, e.g., Section 5 of \cite{BHW}):
Given a graph $G$, let $G'$ be the bipartite graph with vertex set $\V{G'} = \V{G}\times \{0,1\}$, where $(u,0)$ and $(v,1)$ are adjacent in $G'$ whenever either $(u,v)\in \E{G}$ or $u=v$. That is, $G'$ is $G\times K_2$ with the extra edges $((u,0),(u,1))$ for all $u\in \V{G}$. We will show that
$$\hom(G,\HWR)=\hom(G',\Hind).$$
Indeed, consider an independent set $I$ in $G'$. Color $u\in \V{G}$ blue if $(u,1)\in I$, red if $(u,0)\in I$, and white if it is neither red or blue. Note that since $I$ was an independent set and $((u,0),(u,1))\in \E{G'}$, the color of vertex $u$ is well-defined and this coloring is in fact a Widom-Rowlinson coloring of $G$. This same construction also works in the other direction, so
$$\hom(G,\HWR)=\hom(G',\Hind).$$
If $G$ is $d$-regular then $G'$ is $(d+1)$-regular, and $K_{d+1}' = K_{d+1,d+1}$. Applying J. Kahn's result \cite{Kahn} for $(d+1)$-regular bipartite graphs, we see that if $G$ has $n$ vertices then
\begin{align*}
  \hom(G,\HWR)&=\hom(G',\Hind)\\
  &\leq \hom(K_{d+1,d+1},\Hind)^{2n/(2(d+1))}=\hom(K_{d+1},\HWR)^{n/(d+1)}.
\end{align*}
We remark that the transformation $G\to G'$ is also mentioned in \cite{KLL}.

\section{Extension} \label{extension}

In this section we would like to point out that for every graph $H$ there is an $\tilde{H}$ such that
$$\hom(G,\tilde{H})=\hom(G',H),$$
where $G'$ is the bipartite graph defined in the previous section. Exactly the same argument we used for $\HWR$ will work for any graph $\tilde{H}$ constructed in this manner. Actually, the situation is even better. To give the most general version we need a definition.

\begin{Def} Let $G$ be a bipartite graph. Let $H$ be another bipartite graph equipped with a weight function $\nu: \V{H}\to \mathbb{R}_+$. Let $\mathbb{I}_{\E{H}}: \A{H}\times \B{H}\to \{0,1\}$ denote the characteristic function of $\E{H}$. Define 
$$Z_b(G,H)=\sum_{\substack{\varphi:\V{G}\to \V{H}\\\varphi(\A{G})\subseteq \A{H}\\\varphi(\B{G})\subseteq \B{H}}}\prod_{(a,b)\in \E{G}} \mathbb{I}_{\E{H}}(\varphi(a),\varphi(b))\prod_{w \in \V{G}}\nu(\varphi(w)),$$
(The subscript $b$ stands for bipartite.)
If $G$ and $H$ are not necessarily  bipartite graphs, but $H$ is a weighted graph we can still define 
$$Z(G,H)=\sum_{\varphi: \V{G}\to \V{H}}\prod_{(u,v)\in \E{G}}\mathbb{I}_{\E{H}}(\varphi(u),\varphi(v))\prod_{w \in \V{G}}\nu(\varphi(w)).$$
\end{Def}
In the language of statistical phsyics, $Z_b(G,H)$ and $Z(G,H)$ are \emph{partition functions}. 

Somewhat surprisingly, J.~Kahn's result holds even in this general case, as shown by D.~Galvin and P.~Tetali \cite{GalTet}.

\begin{Th} \label{gen3} \cite{GalTet} For any bipartite graph $H$ equipped with the weight function $\nu: \V{H}\to \mathbb{R}_+$ and $\mathbb{I}_{\E{H}}: \A{H}\times \B{H}\to \{0,1\}$, and for any $d$-regular simple graph $G$ on $n$ vertices,
$$Z_b(G,H)\leq Z_b(K_{d,d},H)^{n/(2d)}.$$
\end{Th}

The key observation is that for a bipartite graph $H$ equipped with the weight function $\nu: \V{H}\to \mathbb{R}_+$ and characteristic function $\mathbb{I}_{\E{H}}: \A{H}\times \B{H}\to \{0,1\}$, we can define a weighted graph $\tilde{H}$ with weight function $\tilde{\nu}$ and characteristic function $\mathbb{I}_{\E{\tilde{H}}}$ such that
\begin{align}\label{eqn:transformation}
  Z(G,\tilde{H})=Z_b(G',H)\,,
\end{align}
for any graph $G$ (where $G'$ is the modification of $G$ defined in the previous section).
Indeed, construct $\tilde{H}$ with vertex set $\A{H}\times \B{H}$, edges
$$\mathbb{I}_{\E{\tilde{H}}}((a_1,b_1),(a_2,b_2))=\mathbb{I}_{\E{H}}(a_1,b_2)\mathbb{I}_{\E{H}}(a_2,b_1),$$
and weight function
$$\tilde{\nu} (a,b)=\nu(a)\nu(b)\mathbb{I}_{\E{H}}(a,b).$$
In effect, the vertex set of $\tilde{H}$ is only the edges of $H$ (since non-edge pairs are given weight $0$).
Now, for a map $\varphi:G'\to H$, we can consider the map $\tilde{\varphi}: G\to \tilde{H}$ given by
$$\tilde{\varphi}(u)=(\varphi((u,0)),\varphi((u,1))).$$
By the construction of the graphs $G'$ and $\tilde{H}$, the contribution of $\varphi$ to $Z_b(G,H)$ is the same as the contribution of $\tilde{\varphi}$ to $Z(G,\tilde{H})$, and the result \eqref{eqn:transformation} follows.

Finally, applying Theorem~\ref{gen3} to the $(d+1)$-regular graph $G'$ yields
$$Z(G,\tilde{H})=Z_b(G',H)\leq Z_b(K_{d,d},H)^{2n/(2(d+1))}=Z(K_{d+1},\tilde{H})^{n/(d+1)}.$$

Hence we have proved the following theorem.

\begin{Th} \label{weighted-main} For a bipartite graph $H=(A,B,E)$ with vertex weight function $\nu: \V{H}\to \mathbb{R}_+$ let 
$\tilde{H}$ be the following weighted graph: its vertex set is $E(H)$, its edge set is defined by
$((a_1,b_1),(a_2,b_2))\in E(\tilde{H})$ if and only if $(a_1,b_2)\in E(H)$ and $(a_2,b_1) \in E(H)$,
and the weight function on the vertex set is $\tilde{\nu} (a,b)=\nu(a)\nu(b)$ for $(a,b)\in E(H)$.
Then for any $d$--regular simple graph $G$ on $n$ vertices we have
$$Z(G,\tilde{H})\leq Z(K_{d+1},\tilde{H})^{n/(d+1)}.$$
\end{Th}

We can obtain Conjecture 3 of \cite{CPT} as a corollary by applying this theorem in the case where $H$ is the path on $4$ vertices, $a_1 b_1 a_2 b_2$, with appropriate vertex weights. Indeed, if $\nu(a_1)=1$, $\nu(b_1)=\lambda_b$, $\nu(a_2)=\frac{\lambda_w}{\lambda_b}$, $\nu(b_2)=\frac{\lambda_r\lambda_b}{\lambda_w}$ then $\tilde{H}$ is precisely the Widom-Rowlinson graph with vertex weights $\lambda_b,\lambda_r,\lambda_w$. This proves that even for the vertex-weighted Widom-Rowlinson graph we have
$$Z(G,\HWR)\leq Z(K_{d+1},\HWR)^{n/(d+1)}.$$
Hence we have proved the following theorem.

\begin{Th} Let $\HWR$ be the Widom-Rowlinson graph with vertex weights $\lambda_b,\lambda_w,\lambda_r$. Then 
for any $d$--regular simple graph $G$ on $n$ vertices we have
$$Z(G,\HWR)\leq Z(K_{d+1},\HWR)^{n/(d+1)}.$$
\end{Th}

Now let us consider the special case when $H$ is unweighted ($\nu\equiv 1$). In this case $\tilde{\nu}$ is just $\mathbb{I}_{\E{H}}$, so we can think of $\tilde{H}$ as an unweighted graph with vertex set $\V{\tilde{H}} = \E{H}$. There is an edge in $\tilde{H}$ between edges $e=(a_1, b_1)$ and $f=(a_2,b_2)$ of $H$ whenever $(a_1,b_2)$ and $(a_2,b_1)$ are both also edges of $H$.  This is always the case when either $a_1=a_2$ or $b_1 = b_2$, so in particular every edge $e\in \E{H}=\V{\tilde{H}}$ has a self-loop in $\tilde{H}$, and every pair of incident edges in $H$ are adjacent in $\tilde{H}$. We also get an edge $(e,f)\in \E{\tilde{H}}$ if four vertices $a_1 b_1 a_2 b_2$ are all distinct and form a 4-cycle with $e$ and $f$ as opposite edges. In other words, $\tilde{H}$ is precisely the extended line graph of $H$. Hence as a corollary of Theorem~\ref{weighted-main} we have proved Theorem~\ref{gen}.

If $H$ does not contain any $4$-cycle, then $\tilde{H}$ is simply the line graph of $H$ with loops at every vertex.
In particular, if $H$ is a path (or even cycle of length at least $6$) then $\tilde{H}$ is again a path (or even cycle of length at least $6$), but now with a loop at every vertex. Letting $H^o$ denote the graph obtained by adding a loop at every vertex of the graph $H$, we can write the corollary
\begin{Cor}
  If $H = C_k^o$ (for $k\geq 6$ even) or if $H = P_k^o$ (for any $k$), then for any $d$-regular graph $G$
  $$p_H(G) \leq p_H(K_{d+1}).$$
\end{Cor}

It is a good question how to characterize all of the graphs $\tilde{H}$ which can be obtained this way. Note that since $\tilde{H}$ is always fully-looped, this class has no intersection with the class of graphs found by Galvin \cite{Gal1}: the set of graphs $H_q^{\ell}$ obtained from a complete looped graph on $q$ vertices with $\ell \geq 1$ loops deleted.

\begin{Rem} Let $S_k$ be the star on $k$ vertices. One can show (for details see \cite{Gal1}) that, for large enough $d$,
$$p_{S_k^o}(K_{d+1})<p_{S_k^o}(K_{d,d})$$
for $k\geq 6$. From this example we can see that in order to have $p_H(G) \leq p_H(K_{d+1})$ it is not sufficient merely for $H$ to have a loop at every vertex.
\end{Rem}

L.~Sernau \cite{S} introduced many ideas for extending certain inequalities to a larger class of graphs. For instance, recall that the $H_1\times H_2$ has $\V{H_1\times H_2} = \V{H_1}\times \V{H_2}$ and $((a_1,b_1),(a_2,b_2))\in \E{H_1\times H_2}$ if and only if $(a_1,a_2)\in \E{H_1}$ and $(b_1,b_2)\in \E{H_2}$. Sernau noted that if $H_1$ and $H_2$ are graphs such that
$$p_{H_i}(G)\leq p_{H_i}(K_{d+1})\,,$$
for $i=1,2$, then it is also true that
$$p_{H_1\times H_2}(G)\leq p_{H_1\times H_2}(K_{d+1}).$$
This inequality simply follows from the identity
$$\hom(G,H_1\times H_2) = \hom(G,H_1)\hom(G,H_2),$$
which is explained in \cite{S}. Surprisingly, this observation does not allow us to extend our result to any new graphs, because the product of two extended line graphs is again an extended line graph:
$$\tilde{H}_1\times \tilde{H}_2=\tilde{H}_{12},$$
where $H_{12}=(\A{H_1}\times \A{H_2},\B{H_1}\times \B{H_2},\E{H_1}\times \E{H_2})$.

\section{On a theorem of L.~Sernau}

Theorem 3 of \cite{S} also provides a class of graphs for which $K_{d+1}$ is the maximizing graph. Below we explain the relationships between our results and his theorem.

\begin{Def} Let $H$ and $A$ be graphs. Then the graph $H^A$ is defined as follows: its vertices are the maps $f:V(A)\to V(H)$ and the  $(f_1,f_2)\in E(H^A)$ if $(f_1(u),f_2(v))\in E(H)$ whenever $(u,v)\in E(A)$.
\end{Def}

Then Sernau proved the following theorem.

\begin{Th}\cite{S} Let $G$ be a $d$--regular graph, and let $F=l(H^B)$, where $H$ is an arbitrary graph, $B$ is a bipartite graph, and $l(H^B)$ is the graph induced by the vertices of $H^B$ which have a loop. 
Then
$$p_F(G)\leq p_F(K_{d+1}).$$
\end{Th}

When $H=H_{ind},B=K_2$ then $l(H^B)=H_{WR}$ so this also proves the conjecture of D. Galvin. Note that when $B=K_2$ then $l(H^B)$ is the extended line graph of $H\times K_2$. It is not a great surprise that these results are similar, even the proofs behind these results are strongly related to each other. 

\section{Conjectures}

Let $H$ be a simple graph, i.e., with no multiple edges or loops. Let $H^o$ denote the graph obtained by adding a loop at each vertex of $H$ (so for instance $C_n^o$ denotes the $n$-cycle with a loop at each vertex). 

\begin{Conj} Let $G$ be  a $d$-regular simple graph. Then for any $n\geq 4$
$$p_{C_n^o}(G)\leq p_{C_n^o}(K_{d+1}).$$
\end{Conj}

\begin{Conj} Let $G$ be  a $d$-regular simple graph. Then for any $d\geq 4$
$$p_{S_4^o}(G)\leq p_{S_4^o}(K_{d+1}).$$
Furthermore, for $k\geq 6$
$$p_{S_k^o}(G)\leq p_{S_k^o}(K_{d,d}).$$
\end{Conj}

Finally, for an arbitrary graph $H$ it is not clear how to characterize the maximizers over all $d$-regular graphs $G$ of $p_H(G)$.  If we restrict to bipartite $G$, however, D.~Galvin and P.~Tetali proved that $p_H(G) \le p_H(K_{d,d})$ \cite{GalTet}.  We conjecture that this can be extended to the class of triangle-free graphs.

\begin{Conj} Let $G$ be a $d$--regular triangle-free graph. Then for any graph $H$ we have
$$p_H(G)\leq p_H(K_{d,d}).$$

\end{Conj}

\noindent \textbf{Acknowledgments.} We thank David Galvin and Luke Sernau for helpful conversations. We are also grateful to the anonymous referees for their careful reading and useful suggestions on the paper.

\end{document}